\documentclass[a4paper,10pt]{amsart}

\usepackage{amsmath}
\usepackage{amssymb}
\usepackage{graphicx}

\usepackage{amsthm}
\theoremstyle{definition}
\newtheorem*{dfn1}{Definition     1}

\newtheorem*{exa}{Counterexample }

\newtheorem*{prop1}{Proposition  1}
\newtheorem*{prop2}{Proposition  2}
\newtheorem*{prop3}{Proposition   3}
\newtheorem*{prop4}{Proposition   4}

\newtheorem*{lem1}{Lemma      1}
\newtheorem*{lem2} {Lemma      2}
\newtheorem*{lem3} {Lemma       3}

\newtheorem*{remark}{Remark  }

\newtheorem*{question1}{Question   1}
\newtheorem*{question2}{Question   2}
\newtheorem*{question3}{Question   3}
\newtheorem*{question4}{Question   4}
\newtheorem*{question5} {Question  5}
\newtheorem*{question6} {Question  6}
\newtheorem*{question7} {Question  7}
\newtheorem*{question8} {Question  8}

\theoremstyle{plain}

\newtheorem*{coro1}{Corallary  1}
\newtheorem*{coro2}{Corallary 2 }
\newtheorem*{coro3}{Corollary  3}

\begin{document}

\title[Note On  $Z^{*}$  Algebras]{A Note  On  $Z^{*}$  Algebras }
\author{Ali Taghavi}

\address{Faculty of Mathematics and Computer Science,  Damghan  University,  Damghan,  Iran.}
\email{taghavi@du.ac.ir}

\date{\today}


\begin{abstract}
We  study some  properties of $Z^{*}$  algebras, those  $C^{*}$  algebras which  all positive  elements are zero  divisors.  By  means of an example we show  that an  extension of  a $Z^{*}$ algebra by  a  $Z^{*}$  algebra is not    necessarily  a $Z^{*}$  algebra. However we prove that the  extension of a non $Z^{*}$  algebra  by  a  non $Z^{*}$  algebra is  a non $Z^{*}$  algebra. We also prove  that the tensor product of  a $Z^{*}$  algebra by a  $C^{*}$ algebra is a $Z^{*}$  algebra.\\
As an indirect  consequence of   our  methods  we prove the following inequality type results:\\
i)Let $a_{n}$  be  a  sequence of  positive elements of a $C^{*}$ algebra $A$ which converges  to zero. Then there are positive sequences $b_{n}$ of real numbers and $c_{n}$ of  elements of $A$ which converge to zero such that $a_{n+k}\leq b_{n}c_{k}$.\\
ii)Every   compact   subset of  the  positive  cones of  a  $C^{*}$  algebra has  an upper bound in the algebra.
\end{abstract}

\maketitle
\section*{$Z^{*}$  algebras}
  Throughout  the paper a zero divisor in  a  $C^{*}$  algebra is meant  a right or left zero divisor. However every positive zero divisor is automatically a two sided zero divisor. It is  well known that every  element of  a non unital $C^{*}$ algebra is a topological zero divisor, see \cite{BD}. But it is not true that every non unital  $C^{*}$  algebra  satisfies in the property that all its elements are zero divisor. In this  paper we are interested in $C^{*}$ algebras with this property.\\
  A  $Z^{*}$  algebra  is  a  $C^{*}$  algebra which all positive  elements are zero divisor.
 We will prove  that, for  a  $C^{*}$ algebra $A$, this property  is  equivalent to say that for every compact subset $K$ of of  $A$, there exist  a positive element $x$ with $xK=Kx=0$. In particular,  every element of  a  $Z^{*}$  algebra is  a two sided zero divisor.\\
  The algebra of all compact operators on a non separable  Hilbert space is  an example of  a simple  $Z^{*}$  algebra. Because for a non separable  Hilbert space $H$ and a positive operator $T$ in $K(H)$,   $T$ has  a non trivial kernel. Then for every finite dimensional projection $P$ on a subspace of kernel of $T$ we have $PT=TP=0$. The  algebra of all continuous  functions vanishing at infinity on the long line is  another  example of a $Z^{*}$  algebra, which is also a projection less algebra.
 Direct sum of an uncountable  family of $C^{*}$ algebra is  again  a $Z^{*}$ algebra.\\
The  following lemma  is frequently used in the paper. Its  proof is  a  consequence of lemma I.5.2 in \cite{DAVID} which says that if $0\leq a \leq b$ then  $\sqrt[4]{b}$ is  a factor of $\sqrt{a}$.
  \begin{lem1}
 Let $a,b,c$ be  positive elements of  a $C^{*}$  algebra and $a\leq b$  with $bc=0$. Then $ac=0$.
 \end{lem1}
 We give a  commutative  interpretation for $Z^{*}$ algebras. Recall that every commutative $C^{*}$  algebra is in the  form $C_{0}(X)$, the  algebra of  all continuous functions on a locally compact  Hausdorff space $X$ which vanish at infinity. For the  commutative  interpretation we need to the  following  definition:

\begin{dfn1}
A topological  space $X$  is  approximately $\sigma$ compact, briefly $A$\begin{Large}${\sigma}$\end{Large}$C$, if  it has a sequence of  compact  subsets with dense union.
\end{dfn1}
Now we have the  following  characterization of  all  commutative $Z^{*}$  algebras:

\begin{prop1}
The  commutative  $C^{*}$  algebra $C_{0}(X)$ is  a  $Z^{*}$  algebra if  and  only  if $X$ is  not  an $A$\begin{Large}${\sigma}$\end{Large}$C$  space.
\end{prop1}
\begin{proof}
Let $X$  be an $A$\begin{Large}${\sigma}$\end{Large}$C$  space. Then there are compact  subsets  $K_{n}$  such that $\bigcup_{n=1}^{\infty} K_{n}$ is  dense in $X$. Let $\widetilde{X}=X \cup \{\infty\}$ be the one  point  compactification  of $X$. By urison lemma there is a  sequence of continuous functions $f_{n}:\widetilde{X}\rightarrow [0, 1/n^{2}]$  such that
$f_{n}(K_{n})=1/n^{2}$ and $f_{n}(\infty)=0$. Let $g$ be the restriction of  $\sum_{n=1}^{\infty} f_{n}$ to $X$. Then  $g\in C_{0}(X)$ and $g$  does not  vanish on a  dense subset $\cup K_{n}$. This  obviously  shows that $g$ is  not  a zero divisor hence $C_{0}(X)$ is  not a $Z^{*}$  algebra.  Because the  interior of $h^{-1} \{0\}$  is   not  empty for  every  zero  divisor  $h\in C_{0}(X)$.  For proof  of  the  converse, assume that $g\in C_{0}(X)$ is  not a zero divisor so the interior of $g^{-1}\{0\}$ is empty. Define  a sequence of compact sets  $K_{n}=g^{-1} \{z\in \mathbb{C} \mid 1/n\leq \parallel z \parallel \leq n \}$.
Then each $K_{n}$ is  a compact set and $\cup K_{n}$is dense in $X$. So $X$ is  a $A$\begin{Large}${\sigma}$\end{Large}$C$  space.
\end{proof}

An immediate  consequence of the definition of an $A$\begin{Large}${\sigma}$\end{Large}$C$  space is that every open and  dense subset of a non $A$\begin{Large}${\sigma}$\end{Large}$C$  space is a  non $A$\begin{Large}${\sigma}$\end{Large}$C$  space. Now in  proposition 2 below  we  prove  a non commutative  analogy for this  statement. Recall that a  closed two  sided ideal $J$  of  a  $C^{*}$  algebra  is essential if $J \cap I\neq \{0\}$ for  every  non zero ideal $I$.  Essential ideals  are the  non commutative  analogy of open and dense subsets. The  following  lemma  gives an equivalent condition for an ideal to be  essential:
\begin{lem2}
An ideal $J$ in  a  $C^{*}$ algebra $A$  is  essential if and only if for every positive  element $b\in A$ there exist   $x\in J$ with $x\in (0, b)$.
\end{lem2}

\begin{proof}
Let $J$ be  an essential ideal and $b$ is a  positive element of $A$ with unique positive  square root $c$. Since $J$ is  an  essential ideal there is  a  positive  element $y\in J$ with $\parallel y\parallel \leq 1$ such that $cyc\neq 0$. Then $x=cyc$ is  the  desired element as in the lemma.
Now assume that  $J$ is  an ideal which  satisfies in the  condition of the lemma. Let $I$  be  a non zero ideal in $A$.
Choose  a  positive  element $b\in I$ and  assume  that $x\in J$ belongs to $(0, b)$  then $x$ is  a  non zero element of $I\cap J$, since $I$, being  an ideal of  $A$, is  a hereditary subalgebra of $A$. So $0<x<b \in I$ implies that $x\in I$.

\end{proof}
\begin{prop2}
Every  essential ideal of  a  $Z^{*}$  algebra is  a  $Z^{*}$ algebra.
\end{prop2}

\begin{proof}
Let $J$ be  an essential ideal in a  $Z^{*}$  algebra $A$. Choose  a positive $a\in J$. Then there is  a positive $b\in A$ with $ab=0$. By the  above  lemma there exist $x\in J \cap (0, b)$. So $ax=0$ by lemma 1. This  shows that $J$ is  a  $Z^{*}$  algebra.
\end{proof}

Let $X$  and $Y$ be  two  locally compact  Hausdorff   spaces. Obviously  if either $X$  or $Y$ is  not  an $A$\begin{Large}${\sigma}$\end{Large}$C$  space then $X\times Y$ is  not an $A$\begin{Large}${\sigma}$\end{Large}$C$ space, too. On the other hand the operator theoretical  analogy of the product topology is the spatial or minimal tensor product. This  situation is  a motivation for the  following proposition:
\begin{prop3}
The tensor product of a $Z^{*}$ algebra  by  a  $C^{*}$  algebra is  a  $Z^{*}$  algebra?
\end{prop3}
 The proposition is  an immediate consequence of lemma 1  and lemma 3 below. We  thank Narutaka  Ozawa, for his  proof of lemma 3  in \cite{OZ}.
\begin{lem3}
Let $x$  be  a positive element of $A\otimes_{\text{min}} B$. Then there is  a  single tensor $a\otimes b$ with
$x\leq a\otimes b$.
\end{lem3}
\begin{proof}
In this proof we do not  assume  that the underline  algebras are unital but the  computation is placed in the  corresponding unitization. Let $C$  be  a  $C^{*}$ algebra  and $\sqrt{h_{\alpha}}$  is  an approximate  identity for $C$ with $0<h_{\alpha} \leq 1$. For  a self adjoint  element $y\in C$, define $s(y)=$ the  maximum value of spectrum of $y$. Then for every real number $\lambda \geq s(y)$ and for ever $\epsilon > 0$, there is an $h_{\alpha}$ with
$s(y - \lambda h_{\alpha}) \leq \epsilon$. Because $s(y) \leq \lambda$ implies that $y\leq \lambda$ hence
$\sqrt{h_{\alpha}} y \sqrt{h_{\alpha}} \leq \lambda h_{\alpha}$. So $y-\lambda h_{\alpha} \leq y-\sqrt{h_{\alpha}} y \sqrt{h_{\alpha}}$ therefore  $s(y-\lambda h_{\alpha}) \leq s( y-\sqrt{h_{\alpha}} y \sqrt{h_{\alpha}})$. Since the  net of  real numbers $s( y-\sqrt{h_{\alpha}} y \sqrt{h_{\alpha}})$ approach to zero, we conclude that there is  an $h_{\alpha}$  with $s(y - \lambda h_{\alpha}) \leq \epsilon$. We  apply the  above  statements  to $A\otimes_{\text{min}} B$. Assume that $\sqrt{e_{\alpha}}$, $\sqrt{f_{\alpha}}$ are approximate identity for $A$  and $B$ respectively  with $0\leq e_{\alpha} \leq 1$ and $0\leq f_{\alpha} \leq 1$. So for $h_{\alpha}=e_{\alpha} \otimes f_{\alpha}$, $\sqrt{h_{\alpha}}$ is  an approximate identity for $A\otimes_{\text{min}} B$. let $x$  be  a positive  element in $A\otimes_{\text{min}} B$. Without lose of  generality we can assume that $x\leq 1$. Define  a  sequence of self adjoint  elements $y_{n}$ recursively,  as follows: $y_{0}=x$,\; $y_{n+1}=y_{n}-4^{-n}(e_{n} \otimes f_{n})$ with $s(y_{n}) \leq 4^{-n}$, where $e_{n}$  and $f_{n}$ are chosen from the above  $e_{\alpha}, f_{\alpha}\;'s$ correspond to $\epsilon =4^{-n}$. On the other hand $s(y_{n+1}) \leq 4^{-(n+1)}$ implies that $y_{n+1} \leq 4^{-(n+1)}$. Then for all $n\in \mathbb{N}$ we have  $x=y_{0} \leq  4^{-(n+1)}+\sum_{k=0}^{n} 4^{-k} e_{k} \otimes f_{k} $. So $x\leq e\otimes f$ where $e=\sum 2^{-n}e_{n}$
and $f=\sum 2^{-n}f_{n}$.
\end{proof}

The  following   corollaries  are a consequence of the above lemma:
 \begin{coro1}
 Let $K$  be  a  compact  subset of positive cones of  a  $C^{*}$  algebra $B$. Then $K$ has  an upper bound in $B$.
 \end{coro1}
 \begin{proof}
 Put $A=C(K)$. Let $x:K\rightarrow B$ be  the inclusion map.Then  $x$ is  a positive element of $C(K,B)\simeq A\otimes B$.
 By the  above  lemma, there is  a positive $f\in C(K)$ and  $b\in B$  such that $x(t)\leq f(t)b,\;\text{for all t}\in K$.
 Then $mb$ is  an upper bound  for $K$, where $m=\sup f(t),\;t\in K$
 \end{proof}
 \begin{coro2}
 Assume  that $K$ is  a compact subset   of  a  $Z^{*}$  algebra $A$. Then there  is  a positive  element $x\in A$ such that $xK=Kx=0$.
 \end{coro2}
 \begin{proof}
 Put $\widetilde{K}=\{y^{*}y+ yy^{*} \mid y\in K\}$. Then $\widetilde{K}$ is  a compact  subset of the positive  cones of $A$.
 The  above corollary implies   that there is  a positive   $b\in A$  which is  an upper bound  for $\widetilde{K}$. Since  $A$ is  a $Z^{*}$ algebra, there is  a positive   $x\in A$ with $xb=bx=0$. Using lemma 1 we  conclude  that $x\widetilde{K}=\widetilde{K}x=0$. Hence $xK=Kx=0$, because $y^{*}yx+yy^{*}x=0$ implies that $yx=xy=0$, by lemma 1.

 \end{proof}
 \begin{coro3}
 Let $a_{n}$  be  a  sequence of  positive elements of a $C^{*}$ algebra $A$ which converges to zero. Then there are positive sequences $b_{n}$ of real numbers and $c_{n}$ in $A$ which converge to zero and $a_{n+k}\leq b_{n}c_{k}$
 \end{coro3}
 \begin{proof}
 Put $A=C_{0}(\mathbb{N}),\; B=c_{0}(A)$. Define the  map $f\in C_{0}(\mathbb{N},B)\simeq  C_{0}(\mathbb{N})\otimes B$  with
 \begin{center}
  $f(n)=(a_{n},a_{n+1}, \ldots ) $
  \end{center}
By lemma 3 there are positive  elements $g\in C_{0}(\mathbb{N})$  and  $h\in B$ such that\\
 $f(n)\leq g(n)h$. Let $h$ be  in the  form $h=(c_{1},c_{2},\ldots)$ and put $b_{n}=f(n)$. Then $a_{n+k}\leq b_{n}c_{k}$.
  \end{proof}
 \section*{$Z^{*}$ algebras  and extension theory}
 Let $P$ be  a property about  $C^{*}$ algebras .  A  $C^{*}$  algebra with this  property is  called a
 $P\_ \;\text{algebra}$. A natural question in extension theory is that whether an extension of  a
 $P\_ \;\text{algebra}$ by  a  $P\_ \;\text{algebra}$ is  again a  $P\_ \;\text{algebra}$? Regarding this question, there are some "examples"  and  "non examples". For instance each of the properties of being commutative, unital,  approximately finite dimension,  nuclear or amenable is  an "example" of this situation. But the property of being real rank zero is a "non example". \\
 From this point of view, we consider the property of being $Z^{*}$ algebra. Assume that $B$ and $C$ are two $Z^{*}$ algebras and we have  an extension of $C^{*}$ algebras as follows:\\
 \begin{center}
 $0\rightarrow B\rightarrow A\rightarrow C \rightarrow 0$
 \end{center}
 Does  this  situation implies that $A$ is  a  $Z^{*}$  algebra, too? Equivalently, we have  a $C^{*}$  algebra $A$ and $J$ is a closed two sided ideal in $A$. Assume that $J$ and $A/J$  are $Z^{*}$  algebras. Does this implies that $A$ is  a
 $Z^{*}$  algebra? With the following  formal algebraic computation, we expect that the answer is affirmation: For  convenience assume that $A$ is  commutative and $a\in A$ is  a positive  element. We have to find a non zero positive $y\in A$ with $ay=0$. Since  $A/J$ is  a  $Z^{*}$  algebra, there is  a positive   $b\in A$ such that $ab \in J$. Now $J$ is  a $Z^{*}$ algebra then there is  a positive $c\in J$  such that $abc=0$. Then for $y=bc$ we have $ay=0$. But the problem is that we are not sure that $y$ is  a non zero element. For a given $a$, it is possible  that for all $b$ with $ab\in J$  and  all $c\in J$  with $abc=0$ we necessarily have $bc=0$. The  following example shows that this problem can be hold in certain commutative $C^{*}$algebras:
 \begin{exa}
We give  an example  of a commutative non $Z^{*}$  algebra $A$  which  has an ideal $J$  such that $J$  and $A/J$ are  $Z^{*}$  algebras. Our example is $A=C_{0}(X)$, where $X$ is   the deleted Tychonoff plank, see example 86 in  \cite{COUNTER}. Let $\Omega$  be the first uncountable ordinal  and $\omega$  be the first countable  ordinal. $[0,\;\Omega]$ and $[0,\;\omega]$, as ordered topological spaces,  are compact  spaces. The  deleted Tychonoff plank is $X=[0,\;\Omega]\times [0,\;\omega]-\{(\Omega,\;\omega)\}$. $[0,\;\omega]$ is in the  form $\{0,1,2,\ldots \infty=\omega \}$. We  show that $X$ is  an approximately  $\sigma$ compact  space. For every finite $n\in \mathbb{N}$, define  $K_{n}=[0,\;\Omega]\times \{n\}$. Each $K_{n}$ is  a compact set and $\cup_{n=1}^{\infty} K_{n}$ is  a dense subset of $X$. So $X$ is an $A$\begin{Large}${\sigma}$\end{Large}$C$  space.
Now define an open set $U\subseteq X$  with $U=[0,\;\Omega)\times [0,\;\omega)$. Put $F=X-U$. Then
$[0,\;\Omega) \times \{1\}$ is  a clopen subset of $U$ and $[0,\; \Omega)\times \{ \omega \}$ is  a clopen subset of $F$. On the other hand $[0,\; \Omega)$ is a  non   $A$\begin{Large}${\sigma}$\end{Large}$C$  space and is homomorphic to $[0,\;\Omega) \times \{1\}$ and
$[0,\; \Omega)\times \{ \omega \}$. This  shows that $U$ and $F$ are not $A$\begin{Large}${\sigma}$\end{Large}$C$ spaces. Because every clopen subset of an $A$\begin{Large}${\sigma}$\end{Large}$C$ space is an $A$\begin{Large}${\sigma}$\end{Large}$C$ space. Then $X$  is  an $A$\begin{Large}${\sigma}$\end{Large}$C$ space but $U$  and $X-U$ are not. This  gives us  a  counter example  of   an extension of  a $Z^{*}$  algebra  by  a $Z^{*}$ algebra which resulting extension is not  a $Z^{*}$ algebra.
 \end{exa}
 Despite of this pathology for  $Z^{*}$  algebras, in the  proposition 4 below, we observe that non $Z^{*}$ algebras are well behaved from the extension theoretical  view point. The topological motivation for this  proposition is that if $X$ is  a  locally  compact  space  and $U$ is  an open set in $X$ then $X$ is  an $A$\begin{Large}${\sigma}$\end{Large}$C$   space  provided both $U$  and $X-U$  are $A$\begin{Large}${\sigma}$\end{Large}$C$  spaces:

\begin{prop4}
 An extension of  a non $Z^{*}$ algebra by  a  non  $Z^{*}$ algebra is  again  a non $Z^{*}$  algebra.
 \end{prop4}
\begin{proof}
Let $J$  be  an ideal in $A$ and $\pi: A \rightarrow A/J$ be  the  canonical map. To prove the proposition we equivalently show that if $A$ is  a $Z^{*}$  algebra  but $A/J$ is  not  a  $Z^{*}$  algebra, then $J$ is  a $Z^{*}$  algebra. Let $a\in J$ is  a positive  element. We find  a positive  element $c\in J$  such that $ac=0$. Since $A/J$  is  not  a $Z^{*}$  algebra, there is  a positive  element $b \in A$ such that  $\pi (b)$ is  not  a zero divisor. Then $\pi (a)+ \pi (b)$ is  not  a zero divisor too, by lemma 1. Since $A$ is  a  $Z^{*}$ algebra, there is  a positive element $c \in A$ with  $(a+b)c=c(a+b)=0$. Then $\pi(a+b) \pi (c)=\pi(c) \pi(a+b) =0$. This  shows $\pi(c)=0$, that is
$c\in J$. On the other hand $(a+b)c=c(a+b)=0$ implies  that $ac=ca=0$ by lemma 1. This  completes the  proof  of the proposition.
\end{proof}
\section*{Further questions and remarks}
We present  some  questions and remarks related to the materials of the above  sections: \\

We showed, by  a counter example, that an extension of a $Z^{*}$ algebra by  a $Z^{*}$ algebra is not necessarily a $Z^{*}$  algebra. The underline  algebras in our example were commutative hence non simple. In this line we ask:\\
\begin{question1}
Is an extension of a simple $Z^{*}$ algebra by  a simple $Z^{*}$  algebra,  a $Z^{*}$ algebra?
\end{question1}

In proposition 1 we introduced  a topological interpretation for $Z^{*}$ algebras. Now  we ask a  question about the dynamical interpretation  for  such  algebras:
\begin{question2}
Let $(X,\;f)$  be  a  dynamical system on a  non compact space $X$. This  naturally give's us  an automorphism $\alpha$ on $C_{0}(X)$ and  an action of $\mathbb{Z}$ on $C_{0}(X)$ with $n.g=g\circ f^{n}$. What is  the    dynamical interpretation if  we know that  the crossed product algebra $C_{0}(X) \times_{\mathbb{Z}} \alpha$ is  a  $Z^{*}$  algebra? For  definition of crossed product algebra see \cite[page 222]{DAVID}. Is there a reasonable relation with existence of a  compact invariant  set with non empty interior?
\end{question2}
Our  next question is  about the  representation of $Z^{*}$  algebras. Let $A$  be a  simple  $Z^{*}$ algebra. Then the range of  every representation of $A$ consists only non injective  operators. On the other hand if  a positive element $a$ of a $C^{*}$  algebra is  not  a zero divisor then the linear map on $A$ defined by $x\rightarrowtail ax$ is  an injective operator. We know that such linear maps are the basis of construction of the universal representation of $C^{*}$ algebras on Hilbert spaces. These two situation are  motivations to  ask:
\begin{question3}
Is it true  to say that a simple $C^{*}$  algebra  $A$ is  a  $Z^{*}$  algebra if  and only  if for  every representation
$\pi: A\rightarrow B(H)$, $\pi(A)$ consists only non injective operators?
\end{question3}
\begin{question4}
For  a  $C^{*}$ algebra $A$, what obstructions other than separability , are the most important obstructions for $A$ to have a subalgebra which is  a  $Z^{*}$ algebra? In particular does  $B(H)$, for  a  separable Hilbert space  $H$ contain  a $Z^{*}$ algebra?
\end{question4}
In the next question we are interested to the converse of proposition 3. Assume that $X\times Y$ is not an $A$\begin{Large}${\sigma}$\end{Large}$C$ space. Then either $X$ or $Y$ is necessarily a  non $A$\begin{Large}${\sigma}$\end{Large}$C$ space. So we ask:
\begin{question5}
Assume that $A$ and $B$ are two non $Z^{*}$ algebras. Does it implies that $A \otimes_{\text{min}} B$ is  a non $Z^{*}$ algebra?
\end{question5}

 The proof of proposition 3 was essentially based on lemma 3. So if we could prove a reverse version with the statement  that for every positive element $x\in A \otimes_{\text{min}} B$, there is  a single tensor $a\otimes b$ with $0< a\otimes b \leq x$ then the answer to the above question would be affirmative. But the following example in \cite{Math}  shows that such single positive  tensor is not necessarily exist:\\
 Assume that  $A=B=M_{2}(\mathbb{C})$ so $A \otimes B=M_{4}(\mathbb{C})$.
 Put $x=\left(\begin{smallmatrix}1&0&0&1\\0&0&0&0\\0&0&0&0\\1&0&0&1 \end{smallmatrix}\right)$. It can be shown that there are no non zero matrices $a,b\in  M_{2}(\mathbb{C})$ such that $0< a\otimes b \leq x$.\\
 \begin{question6}
 In this question we search for  a pure algebraic version of $Z^{*}$ algebras.  A commutative ring $R$ is  a $Z\_\; \text{ring}$ if every element of $R$ is  a zero divisor. Can one prove  similar results as propositions 2,3 and 4 in the category of  commutative rings?
 \end{question6}
Note  that every (not necessarily proper)  continuous  image of  an $A$\begin{Large}${\sigma}$\end{Large}$C$ space is  again an $A$\begin{Large}${\sigma}$\end{Large}$C$  space . On the  other  hand the  non commutative  analogues of continuous maps which are not necessarily  proper are the  Woronowicz morphisms. Recall that a Woronowicz morphism from  an algebra $A$ to $B$ is  a morphism $f$ from $A$ to the  multiplier algebra $\mathcal{M}(B)$  such that the  linear span generated by $\{f(a)b \mid a\in A,\;b\in B \}$ is  dense in B, see \cite[page 15]{E}. So it is  natural to ask the  following question:
\begin{question7}
Assume that there is an injective Woronowicz morphism from  a  $Z^{*}$  algebra $A$ to a  $C^{*}$  algebra $B$. Does it implies that $B$ is  a  $Z^{*}$  algebra?
\end{question7}
\begin{question8}
Our  final  question is  about  a generalization of proposition 1 for  commutative Banach  algebras and  non commutative  $C^{*}$  algebras:
\begin{itemize}

\item Let $A$  be  a  commutative  Banach  algebra which all elements  are zero  divisor. Is it true to say the its Gelfand  spectrum   is  not  an $A$\begin{Large}${\sigma}$\end{Large}$C$  space?
\item Let  $A$  be  a non commutative  $Z^{*}$  algebra. What can be  said  about the topology of $\widehat{A}$, the  space of irreducible representations of $A$ or the topology of $P(A)$, the  space  of  pure  states of $A$?
\end{itemize}
\end{question8}
\begin{remark}
Note that the converse of the first part of the  above  question is  not true. There is a  commutative  Banach  algebra $A$ which Gelfand  spectrum is  not an  $A$\begin{Large}${\sigma}$\end{Large}$C$  space but $A$ has no zero divisor. Put $A=\text{the B\_valued disck algebra}$ where $B$ is  a commutative  $Z^{*}$  algebra.   That is the  space of all holomorphic  functions  from the unit disc
$\mathbb{D}$ to $B$ with continuous  extension to $\overline{\mathbb{D}}$. Let $\widehat{B}$ be  the spectrum of $B$.
Then $\widehat{A}=\overline{\mathbb{D}} \times \widehat{B}$ is   not  an $A$\begin{Large}${\sigma}$\end{Large}$C$    space.

\end{remark}

\end{document}